\documentclass{article}%
\usepackage{amsfonts}%
\usepackage{amsmath}%
\usepackage{amssymb}%
\usepackage{graphicx}
\providecommand{\U}[1]{\protect\rule{.1in}{.1in}}
\newtheorem{theorem}{Theorem}

\newtheorem{corollary}[theorem]{Corollary}

\newtheorem{definition}[theorem]{Definition}

\newtheorem{lemma}[theorem]{Lemma}

\newenvironment{proof}[1][Proof]{\noindent\textbf{#1.} }{\ \rule{0.5em}{0.5em}}
\begin{document}

\title{On the wave equation with variable exponent nonlinearity and distributive delay}
\author{Mohammad Kafini\\Department of Mathematics\\The Interdisciplinary Research Center in Construction\\and Building Materials\\KFUPM, Dhahran 31261 \\Saudi Arabia\\E-mail: mkafini@kfupm.edu.sa}
\date{}
\maketitle

\begin{abstract}
In this work, we are concerned with a nonlinear wave equation with variable
exponents. A distributive delay is imposed into the damping term with variable
exponents nonlinearity. Firstly, we show that the global nonexistence time can
be dominated. Secondly, global existence of solutions is shown under some
suitable conditions on the initial data. Finally, the decay rates of that
solutions are established as well.



\end{abstract}

\section{Introduction}

\qquad In this work, we are concerned with the following delayed nonlinear
wave problem
\begin{equation}
\left\{
\begin{tabular}
[c]{lll}%
$u_{tt}-\Delta u+\mu_{1}u_{t}(x,t)|u_{t}|^{m(x)-2}(x,t)$ &  & \\
$+\int_{\tau_{1}}^{\tau_{2}}\mu_{2}\left(  \tau\right)  u_{t}(x,t-\tau
)|u_{t}|^{m(x)-2}(x,t-\tau)d\tau=u\left\vert u\right\vert ^{p(x)-2},$ & in &
$\Omega\times(0,\infty)$\\
$u(x,t)=0,$ & in & $\partial\Omega\times\lbrack0,\infty)$\\
$u(x,0)=u_{0}(x),\emph{\quad}u_{t}(x,0)=u_{1}(x),$ & in & $\Omega$\\
$u_{t}\left(  x,t-\tau\right)  =f_{0}\left(  x,t-\tau\right)  ,$ & in &
$\Omega\times(0,\tau_{2}),$%
\end{tabular}
\ \right. \label{P}%
\end{equation}
\ \ \ \ 

The domain $\Omega$ $\subset\mathbb{R}^{n}$ is bounded with a smooth boundary
$\partial\Omega$. The exponents $m(\cdot)$ and $p(\cdot)$ are given measurable
functions on $\overline{\Omega}$ and satisfy
\begin{equation}
2\leq m_{1}\leq m(x)\leq m_{2}<p_{1}\leq p(x)\leq p_{2}\leq\frac{2\left(
n-1\right)  }{n-2},\text{ \ \ }n\geq3,\label{eq1}%
\end{equation}
where
\begin{align*}
m_{1}  & :=\text{ess inf}_{x\in\Omega}m(x),\text{ \ \ \ \ \ \ \ \ \ \ }%
m_{2}:=\text{ess sup}_{x\in\Omega}m(x),\\
p_{1}  & :=\text{ess inf}_{x\in\Omega}p(x)\text{ \ \ \ and \ \ \ \ }%
p_{2}:=\text{ess sup}_{x\in\Omega}p(x)
\end{align*}
and the log-H\"{o}lder continuity condition:
\begin{equation}
\left\vert q(x)-q(y)\right\vert \leq-\frac{A}{\log\left\vert x-y\right\vert
},\label{eq2}%
\end{equation}
for $x,y\in\Omega,$ with $\left\vert x-y\right\vert <\delta,$ $A>0$ \ and
$0<\delta<1.$

The limits $\tau_{2}>\tau_{1}$ are such that $\mu_{2}:[\tau_{1},\tau
_{2}]\rightarrow\mathbb{R}^{+}$ is a real positive function represents
distributive time delay, satisfies%
\begin{equation}
\int_{\tau_{1}}^{\tau_{2}}\mu_{2}\left(  \tau\right)  d\tau<\mu_{1}%
,\label{eq3}%
\end{equation}
for a positive constant $\mu_{1}.$ The functions $u_{0},$ $u_{1},$ $f_{0}$ are
the initial-boundary conditions and the history data to be specified later.

Problem $\left(  \ref{P}\right)  $ has been extensively studied in the
constant-exponent case. For instance, in the absence of delay ($\mu_{2}=0$),
many blow-up and decay results have been proved (see $\cite{1}-\cite{6}$).
Moreover, the idea of establishing a lower bound for the blow-up time was
presented in several works. See, for instance, $\cite{7}\ $and $\cite{8}$.

In the presence of the delay term ($\mu_{2}\neq0$), the situation is quite
different as it is well known that the delay term can be a source of
instability unless some additional stabilization mechanisms are added. Indeed,
Nicaise and Pignotti $\cite{9}$ examined the wave equation
\begin{equation}
u_{tt}(x,t)-\Delta u(x,t)+a_{0}u_{t}(x,t)+au_{t}(x,t-\tau)=0\text{ \ in
\ }\Omega\times(0,+\infty),\label{eqint1}%
\end{equation}
where $\Omega$ $\subset\mathbb{R}^{n}$ is a bounded domain and $a,$ $a_{0}$
are positive real parameters, and they proved that the system is exponentially
stable under the condition ($0\leq a<a_{0}$). In the case $a\geq a_{0}$, they
produced a sequence of delays for which the corresponding solution of $\left(
\ref{eqint1}\right)  $ is instable. After that, various types of delay were
considered and similar stability results were established. See, in this
regard, $\left(  \cite{10}-\cite{14}\right)  $.

Considering distributive delays, we mention the work of Mustafa \textit{et al}
in $\cite{15}$. The authors considered a thermoelastic system with an internal
distributed delay. They use the energy method to prove, under a suitable
assumption on the weight of the delay, that the damping effect through heat
conduction is strong enough to uniformly stabilize the system even in the
presence of time delay.

Recently, several physical phenomena such as flows of electro-rheological
fluids or fluids with temperature-dependent viscosity, nonlinear
viscoelasticity, filtration processes through a porous media and image
processing are modelled by equations with variable exponents of nonlinearity.
More details on these problems can be found in $\left(  \cite{16}%
-\cite{18}\right)  $.

Regarding hyperbolic problems with nonlinearities of variable-exponent type,
we mention here the work of Antontsev $\cite{19}$, where he considered the
equation
\begin{equation}
u_{tt}-\operatorname{div}\left(  a(x,t)\left\vert \nabla u\right\vert
^{p(x,t)-2}\nabla u\right)  -\alpha\Delta u_{t}=b(x,t)u\left\vert u\right\vert
^{\sigma(x,t)-2},\text{ in }\Omega\times(0,\infty)\label{eqint2}%
\end{equation}
and proved several blow-up results under specific conditions on $a,b,\alpha
,p,\sigma$ and for certain solutions with non-positive initial energy. Also,
in $\cite{20}$, Antontsev proved the existence of local and global weak
solutions of $\left(  \ref{eqint2}\right)  $ by means of Galerkin's
approximations in spaces of Orlicz-Sobolev type. Then he established the blow
up for weak solutions with nonpositive energy functional. Guo and
Gao\ $\cite{21}$ proved blow up for solutions to quasiliear hyperbolic
equations with $p(x,t)$-Laplacian and positive initial energy. Messaoudi and
Talahmeh $\cite{22}$ considered
\begin{equation}
u_{tt}-\operatorname{div}\left(  a(x,t)\left\vert \nabla u\right\vert
^{m(x)-2}\nabla u\right)  -\mu u_{t}=u\left\vert u\right\vert ^{p(x)-2}%
\label{eqint3}%
\end{equation}
and established a blow-up result for solutions to $\left(  \ref{eqint3}%
\right)  $ with an arbitrary positive initial energy. As a matter of fact,
they extended the result established in $\cite{23}$, from the
constant-exponent nonlinearity to variable-exponent nonlinearities. For more
problems involving variable-exponent nonlinearities, we refer the reader to
Antontsev and Shmarev $\cite{24}$, Galaktionov $\cite{25}$, Algarabli
\textit{et al} $\cite{26}$.

In our work, we aim to study the problem $\left(  \ref{P}\right)  $ where the
distributive delay acting in the variable-exponent nonlinearly damping term.
Precisely, we establish a global nonexistence and existence results under
sufficient conditions on $\mu_{1},\mu_{2},m,p$ and the initial data.
Consequantely, we obtain blow-up and decay results. To the best of our
knowledge, there is no work which dealt with the distibutive delay in the
variable-exponent nonlinearity. This paper consists of three sections, in
addition to the introduction. In Section 2, we recall the definitions of the
variable exponent Lebesgue spaces $L^{p(\cdot)}(\Omega),$ the Sobolev spaces
$W^{1,p(\cdot)}(\Omega)$, as well as some of their properties. In section 3,
we transform our problem by introducing a new functional together with its
corresponding energy. In Section 4, we prove a global nonexistence result and
establish a lower bound for the nonexistence time. Section 5 is devoted to the
global existence and the decay result.

\section{Preliminaries}

In this section, we present some materials needed for the statement and the
proof of our results. In what follows, we give definitions and properties
related to Lebesgue and Sobolev spaces with variable exponents.

Let $\Omega$ be a domain of $\mathbb{R}^{n}$ with $n\geq2$ and $p:\Omega
\longrightarrow\lbrack1,\infty)$ be a measurable function (here $p(\cdot) $
has nothing to do with that we used in our problem). The Lebesgue space
$L^{p(\cdot)}(\Omega)$ with a variable exponent $p(\cdot)$ is defined by
\[
L^{p(\cdot)}(\Omega)=\left\{  u:\Omega\longrightarrow\mathbb{R};\text{
measurable in }\Omega\text{ and }\int_{\Omega}\left\vert \lambda
u(x)\right\vert ^{p(x)}dx<+\infty\right\}  ,
\]
for some $\lambda>0.$ The Luxembourg-type norm is given by
\[
\left\Vert u\right\Vert _{p(\cdot)}:=\inf\left\{  \lambda>0:\int_{\Omega
}\left\vert \frac{u(x)}{\lambda}\right\vert ^{p(x)}dx\leq1\right\}  .
\]
The space $L^{p(\cdot)}(\Omega),$ equipped with this norm, is a Banach space
(see $\cite{21}$). The variable-exponent Sobolev space $W^{1,p(\cdot)}%
(\Omega)$ is defined as
\[
W^{1,p(\cdot)}(\Omega)=\left\{  u\in L^{p(\cdot)}(\Omega)\text{ such that
}\nabla u\text{ exists and }\left\vert \nabla u\right\vert \in L^{p(\cdot
)}(\Omega)\right\}  .
\]
This space is a Banach space with respect to the norm
\[
\left\Vert u\right\Vert _{W^{1,p(\cdot)}(\Omega)}=\left\Vert u\right\Vert
_{p(\cdot)}+\left\Vert \nabla u\right\Vert _{p(\cdot)}.
\]
The space $W_{0}^{1,p(\cdot)}(\Omega)$ is defined to be the closure of
$C_{0}^{\infty}(\Omega)$ in $W^{1,p(\cdot)}(\Omega).$ The definition of the
space $W_{0}^{1,p(\cdot)}(\Omega)$ is usually different from the constant
exponent case. However, under condition $\left(  \ref{eq2}\right)  $ both
definitions coincide (see $\cite{27}$). The dual space of $W_{0}^{1,p(\cdot
)}(\Omega)$ is $W_{0}^{-1,p^{\prime}(\cdot)}(\Omega)$ defined in the same way
as in the classical Sobolev spaces, where
\[
\frac{1}{p(\cdot)}+\frac{1}{p^{\prime}(\cdot)}=1.
\]

\begin{lemma}
[\textbf{Poincar\'{e}'s inequality }$\cite{22}$]\emph{Let }$\Omega$\emph{\ be
a bounded domain of }$\mathbb{R}^{n}$\emph{\ and }$p(\cdot)$\emph{\ satisfies
}$\left(  \ref{eq1}\right)  $\emph{, then}
\[
\left\Vert u\right\Vert _{p(\cdot)}\leq C\left\Vert \nabla u\right\Vert
_{p(\cdot)},\text{ for all }u\in W_{0}^{1,p(\cdot)}(\Omega),
\]
\emph{where the positive constant }$C$\emph{\ depends on }$p(\cdot
)$\emph{\ and }$\Omega.$\emph{\ In particular, the space }$W_{0}^{1,p(\cdot
)}(\Omega)$\emph{\ has an equivalent norm given by \qquad}
\[
\left\Vert u\right\Vert _{W_{0}^{1,p(\cdot)}(\Omega)}=\left\Vert \nabla
u\right\Vert _{p(\cdot)}.
\]

\end{lemma}

\begin{lemma}
[$\cite{20}$]\textbf{\ }\emph{If }$p:\overline{\Omega}\longrightarrow
\lbrack1,\infty)$\emph{\ is continuous and }
\[
2\leq p_{1}\leq p(x)\leq p_{2}\leq\frac{2n}{n-2},\text{ \ \ }n\geq3,
\]
\emph{then the embedding }$H_{0}^{1}(\Omega)\hookrightarrow L^{p(\cdot
)}(\Omega)$\emph{\ is continuous.}
\end{lemma}

\begin{lemma}
[$\cite{19}$]\noindent\ \emph{If }$p:\Omega\longrightarrow\lbrack1,\infty
)$\emph{\ is a measurable function and} $p_{2}<\infty,$ \emph{then}
$C_{0}^{\infty}(\Omega)$ \emph{is dense in} $L^{p(\cdot)}(\Omega).$
\end{lemma}

\begin{lemma}
[\textbf{H\"{o}lder's inequality }$\cite{19}$]\noindent\ \emph{Let }%
$p,q,s\geq1$\emph{\ be measurable functions defined on }$\Omega$\emph{\ such
that }
\[
\frac{1}{s(y)}=\frac{1}{p(y)}+\frac{1}{q(y)},\text{ \ \ for a.e. }y\in\Omega.
\]
\emph{If }$f\in L^{p(\cdot)}(\Omega)$\emph{\ and }$g\in L^{q(\cdot)}(\Omega
)$\emph{, then }$fg\in L^{s(\cdot)}(\Omega)$\emph{\ and}%
\[
\left\Vert fg\right\Vert _{s(\cdot)}\leq2\left\Vert f\right\Vert _{p(\cdot
)}\left\Vert g\right\Vert _{q(\cdot)}\emph{.}
\]

\end{lemma}

\begin{lemma}
[\textbf{Unit Ball Property} $\cite{19}$]\emph{Let }$p$\emph{\ be a measurable
function on }$\Omega$\emph{. Then}%
\[
\left\Vert f\right\Vert _{p(\cdot)}\leq1\text{ \ if and only if }%
\varrho_{p(\cdot)}(f)\leq1,
\]
\emph{where}
\[
\varrho_{p(\cdot)}(f)=\int_{\Omega}\left\vert f(x)\right\vert ^{p(x)}dx.
\]

\end{lemma}

\begin{lemma}
[$\cite{20}$]\emph{If }$p$\emph{\ is a measurable function on }$\Omega
$\emph{\ satisfying }$\left(  \ref{eq1}\right)  $\emph{, then}
\[
\min\left\{  \left\Vert u\right\Vert _{p(\cdot)}^{p_{1}},\left\Vert
u\right\Vert _{p(\cdot)}^{p_{2}}\right\}  \leq\varrho_{p(\cdot)}(u)\leq
\max\left\{  \left\Vert u\right\Vert _{p(\cdot)}^{p_{1}},\left\Vert
u\right\Vert _{p(\cdot)}^{p_{2}}\right\}  ,
\]
\emph{for a.e. }$x\in\Omega$\emph{\ and for any }$u\in L^{p(\cdot)}(\Omega).$
\end{lemma}

\section{Energy setting}

In this section, we define the energy functional and show that it is
dissipative. We first introduce the new variable, as in $\cite{9}$,
\[
z(x,\rho,t,\tau)=u_{t}(x,t-\tau\rho),\qquad x\in\Omega,\text{ }\rho
\in(0,1),\text{ }\tau\in(\tau_{1},\tau_{2}),\text{ }t>0.
\]

\noindent Thus, we have
\[
\tau z_{t}(x,\rho,t,\tau)+z_{\rho}(x,\rho,t,\tau)=0,\quad x\in\Omega,\text{
}\rho\in(0,1),\text{ }\tau\in(\tau_{1},\tau_{2}),\text{ }t>0.
\]

\noindent Then, problem $\left(  \ref{P}\right)  \ $takes the form%
\begin{equation}
\left\{
\begin{tabular}
[c]{lll}%
$u_{tt}-\Delta u+\mu_{1}u_{t}(x,t)\left\vert u_{t}(x,t)\right\vert ^{m(x)-2} $
&  & \\
$+\int_{\tau_{1}}^{\tau_{2}}\mu_{2}\left(  \tau\right)  z(x,1,t,\tau
)\left\vert z(x,1,t,\tau)\right\vert ^{m(x)-2}d\tau$ &  & \\
$=u(x,t)\left\vert u(x,t)\right\vert ^{p(x)-2},$ & in & $\Omega\times
(0,\infty)\times(\tau_{1},\tau_{2})$\\
$\tau z_{t}(x,\rho,t,\tau)+z_{\rho}(x,\rho,t,\tau)=0$ & in & $\Omega
\times(0,1)\times(\tau_{1},\tau_{2}),\text{ }t>0.$\\
$z(x,\rho,0,\tau)=f_{0}(x,-\rho\tau),$ & in & $\Omega\times(0,1)\times
(\tau_{1},\tau_{2})$\\
$u(x,t)=0,$ & in & $\partial\Omega\times\lbrack0,\infty)$\\
$u(x,0)=u_{0}(x),\emph{\quad}u_{t}(x,0)=u_{1}(x),$ & in & $\Omega.$%
\end{tabular}
\right. \label{P1}%
\end{equation}

\begin{definition}
\emph{For }$T>0$\emph{\ fixed, we call }$(u,z)$\emph{\ a strong solution if }%
\begin{align*}
u  & \in C^{2}([0,T);L^{2}(\Omega))\cap C^{1}([0,T);H_{0}^{1}(\Omega))\cap
C([0,T);H^{2}(\Omega)\cap H_{0}^{1}(\Omega)),\\
u_{t}  & \in L^{m(\cdot)}(\Omega\times(0,T)),\\
z  & \in C^{1}([0,1]\times\lbrack0,T);L^{2}(\Omega))\cap L^{\infty
}((0,T);L^{m(\cdot)}((0,1)\times\Omega))
\end{align*}
\emph{and satisfies the equations of }$\left(  \ref{P1}\right)  $\emph{\ in
}$H^{-1}(\Omega)$\emph{\ and }$L^{2}(\Omega)$\emph{\ respectively and the
initial data.}
\end{definition}

\noindent See $\cite{11}$ for the well-posedness of a similar problems.

The\ energy functional associated to $\left(  \ref{P1}\right)  $ is given by
\begin{align}
E(t)  & :=\frac{1}{2}||u_{t}||_{2}^{2}+\frac{1}{2}||\nabla u||_{2}^{2}%
+\int_{\Omega}\int_{0}^{1}\int_{\tau_{1}}^{\tau_{2}}\tau\left(  \mu_{2}\left(
\tau\right)  +\xi(x)\right)  \frac{\left\vert z(x,\rho,t,\tau)\right\vert
^{m(x)}}{m(x)}d\tau d\rho dx\nonumber\\
& -\int_{\Omega}\frac{\left\vert u\right\vert ^{p(x)}}{p(x)}dx,\label{eq4}%
\end{align}
for $t\geq0$ and $\xi$ is a positive continuous function satisfying
\begin{equation}
\int_{\tau_{1}}^{\tau_{2}}\mu_{2}\left(  \tau\right)  d\tau+\left(  \tau
_{2}-\tau_{1}\right)  \frac{\xi(x)}{m(x)}<\mu_{1},\text{ \ \ \ }x\in
\overline{\Omega}.\label{eq5}%
\end{equation}
One can take, for instance,
\begin{equation}
\xi(x)=\frac{m(x)}{2\left(  \tau_{2}-\tau_{1}\right)  }\left[  \mu_{1}%
-\int_{\tau_{1}}^{\tau_{2}}\mu_{2}\left(  \tau\right)  d\tau\right]  >0,\text{
\ \ \ on \ \ }\overline{\Omega}.\label{eq6}%
\end{equation}
The following lemma shows that the associated energy of the problem is
nonincreasing under the condition $\left(  \ref{eq5}\right)  $.\noindent

\begin{lemma}
\noindent\emph{Let }$u$\emph{\ be the solution of }$\left(  \ref{P1}\right)  .
$\emph{\ Then, for some }$C_{0}>0$\emph{,}
\begin{equation}
E^{\prime}\left(  t\right)  \leq-C_{0}\left[  \int_{\Omega}\left(  \left\vert
u_{t}\right\vert ^{m(x)}+\int_{\tau_{1}}^{\tau_{2}}\left\vert z(x,1,t,\tau
)\right\vert ^{m(x)}d\tau\right)  dx\right]  \leq0.\label{eq7}%
\end{equation}

\end{lemma}

\noindent\textbf{Proof.} Multiplying equation $\left(  \ref{P1}\right)  _{1}$
by $u_{t}$ and integrating over $\Omega$ and multiplying $\left(
\ref{P1}\right)  _{2}$ by $\left[  \mu_{2}\left(  \tau\right)  +\xi(x)\right]
\left\vert z\right\vert ^{m(x)-2}z$ and integrating over $\Omega
\times(0,1)\times(\tau_{1},\tau_{2})$,\ then summing up, we get
\begin{align}
& \frac{d}{dt}\left[  \frac{1}{2}||u_{t}||_{2}^{2}+\frac{1}{2}||\nabla
u||_{2}^{2}+\int_{\Omega}\int_{0}^{1}\int_{\tau_{1}}^{\tau_{2}}\frac
{\tau\left[  \mu_{2}\left(  \tau\right)  +\xi(x)\right]  \left\vert
z(x,\rho,t,\tau)\right\vert ^{m(x)}}{m(x)}d\tau d\rho dx\right] \nonumber\\
& -\frac{d}{dt}\int_{\Omega}\frac{\left\vert u\right\vert ^{p(x)}}%
{p(x)}dx\label{eq8}\\
& =-\mu_{1}\int_{\Omega}\left\vert u_{t}\right\vert ^{m(x)}dx-\int_{\Omega
}\int_{0}^{1}\int_{\tau_{1}}^{\tau_{2}}\left[  \mu_{2}\left(  \tau\right)
+\xi(x)\right]  \left\vert z\right\vert ^{m(x)-2}zz_{\rho}d\tau d\rho
dx\nonumber\\
& -\int_{\Omega}\int_{\tau_{1}}^{\tau_{2}}u_{t}\mu_{2}\left(  \tau\right)
z(x,1,t,\tau)\left\vert z(x,1,t,\tau)\right\vert ^{m(x)-2}d\tau dx.\nonumber
\end{align}
\noindent We, now, estimate the last two terms of the right-hand side of
$\left(  \ref{eq8}\right)  $ as follows,
\begin{align*}
& -\int_{\Omega}\int_{0}^{1}\int_{\tau_{1}}^{\tau_{2}}\left[  \mu_{2}\left(
\tau\right)  +\xi(x)\right]  \left\vert z(x,\rho,t,\tau)\right\vert
^{m(x)-2}zz_{\rho}(x,\rho,t,\tau)d\tau d\rho dx\\
& =-\int_{\Omega}\int_{0}^{1}\int_{\tau_{1}}^{\tau_{2}}\frac{\partial
}{\partial\rho}\left(  \frac{\left[  \mu_{2}\left(  \tau\right)
+\xi(x)\right]  \left\vert z(x,\rho,t,\tau)\right\vert ^{m(x)}}{m(x)}\right)
d\tau d\rho dx\\
& =\int_{\Omega}\int_{\tau_{1}}^{\tau_{2}}\frac{\left[  \mu_{2}\left(
\tau\right)  +\xi(x)\right]  }{m(x)}\left(  \left\vert z(x,0,t,\tau
)\right\vert ^{m(x)}-\left\vert z(x,1,t,\tau)\right\vert ^{m(x)}\right)  d\tau
dx\\
& =\int_{\Omega}\int_{\tau_{1}}^{\tau_{2}}\frac{\left[  \mu_{2}\left(
\tau\right)  +\xi(x)\right]  }{m(x)}\left\vert u_{t}\right\vert ^{m(x)}d\tau
dx-\int_{\Omega}\int_{\tau_{1}}^{\tau_{2}}\frac{\left[  \mu_{2}\left(
\tau\right)  +\xi(x)\right]  }{m(x)}\left\vert z(x,1,t,\tau)\right\vert
^{m(x)}d\tau dx\\
& =\int_{\Omega}\frac{1}{m(x)}\left(  \int_{\tau_{1}}^{\tau_{2}}\mu_{2}\left(
\tau\right)  d\tau+\xi(x)\left(  \tau_{2}-\tau_{1}\right)  \right)  \left\vert
u_{t}\right\vert ^{m(x)}dx\\
& -\int_{\Omega}\int_{\tau_{1}}^{\tau_{2}}\frac{\left[  \mu_{2}\left(
\tau\right)  +\xi(x)\right]  }{m(x)}\left\vert z(x,1,t,\tau)\right\vert
^{m(x)}d\tau dx.
\end{align*}
For the last term\noindent, we use Young's inequality with $q=\frac
{m(x)}{m(x)-1}$ and $q^{\prime}=m(x)$ to get%
\[
\left\vert u_{t}\right\vert \left\vert z(x,1,t,\tau)\right\vert ^{m(x)-1}%
\leq\frac{m(x)-1}{m(x)}\left\vert u_{t}\right\vert ^{m(x)}+\frac{1}%
{m(x)}\left\vert z(x,1,t,\tau)\right\vert ^{m(x)}.
\]
Consequently, we arrive at%
\begin{align*}
& -\int_{\tau_{1}}^{\tau_{2}}\mu_{2}\left(  \tau\right)  \int_{\Omega}%
u_{t}z\left\vert z(x,1,t,\tau)\right\vert ^{m(x)-2}dxd\tau\\
& \leq\left(  \int_{\tau_{1}}^{\tau_{2}}\mu_{2}\left(  \tau\right)
d\tau\right)  \int_{\Omega}\frac{m(x)-1}{m(x)}\left\vert u_{t}(t)\right\vert
^{m(x)}dx+\int_{\Omega}\int_{\tau_{1}}^{\tau_{2}}\frac{\mu_{2}\left(
\tau\right)  }{m(x)}\left\vert z(x,1,t,\tau)\right\vert ^{m(x)}d\tau dx.
\end{align*}
\noindent Hence, we obtain
\begin{align*}
\frac{dE(t)}{dt}  & \leq-\int_{\Omega}\left[  \mu_{1}-\left(  \int_{\tau_{1}%
}^{\tau_{2}}\mu_{2}\left(  \tau\right)  d\tau+\left(  \tau_{2}-\tau
_{1}\right)  \frac{\xi(x)}{m(x)}\right)  \right]  \left\vert u_{t}%
(t)\right\vert ^{m(x)}dx\\
& -\int_{\Omega}\int_{\tau_{1}}^{\tau_{2}}\frac{\xi\left(  x\right)  }%
{m(x)}\left\vert z(x,1,t,\tau)\right\vert ^{m(x)}d\tau dx.
\end{align*}
\noindent Finally, the relation $\left(  \ref{eq5}\right)  $ yields, $\forall
$\ $x\in\overline{\Omega},$%
\[
f(x)=\mu_{1}-\left(  \int_{\tau_{1}}^{\tau_{2}}\mu_{2}\left(  \tau\right)
d\tau+\left(  \tau_{2}-\tau_{1}\right)  \frac{\xi(x)}{m(x)}\right)  >0.
\]
Since $m(x)$ is bounded, hence $\xi\left(  x\right)  ,$ we deduce that $f(x)$
is bounded. Therefore, for%
\[
C_{0}=\max\left\{  \inf_{\text{\ }x\in\overline{\Omega}}f(x),\inf
_{\text{\ }x\in\overline{\Omega}}\frac{\xi\left(  x\right)  }{m(x)}\right\}
>0,
\]
\noindent we have
\[
E^{\prime}\left(  t\right)  \leq-C_{0}\left[  \int_{\Omega}\left\vert
u_{t}(t)\right\vert ^{m(x)}dx+\int_{\Omega}\int_{\tau_{1}}^{\tau_{2}%
}\left\vert z(x,1,t,\tau)\right\vert ^{m(x)}d\tau dx\right]  \leq
0.\ \ \ \ \ \ \ \ \ \ \
\endproof
\]

\section{Global nonexistence}

In order to prove our global nonexistence result, we assume in addition to
$\left(  \ref{eq4}\right)  $, that $E(0)<0.$ Also, we set
\begin{equation}
H(t)=-E(t),\label{eq9}%
\end{equation}
hence,
\begin{equation}
H^{\prime}(t)=-E^{\prime}(t)\geq0\text{ \ \ and \ \ }0<H(0)\leq H(t)\leq
\int_{\Omega}\frac{\left\vert u\right\vert ^{p(x)}}{p(x)}dx\leq\frac{1}{p_{1}%
}\varrho(u),\label{eq10}%
\end{equation}
where
\[
\varrho(u)=\varrho_{p(\cdot)}(u)=\int_{\Omega}\left\vert u\right\vert
^{p(x)}dx.
\]
We state without proof, the following technical lemmas and corollaries. (See
$\cite{22}$ for proofs)

\begin{lemma}
\noindent\emph{Suppose the condition }$\left(  \ref{eq1}\right)
$\emph{\ holds. Then there exists a positive }$C>1$\emph{, depending on
}$\Omega$\emph{\ only, such that}
\[
\varrho^{\frac{s}{p_{1}}}(u)\leq C\left(  ||\nabla u(t)||_{2}^{2}%
+\varrho(u)\right)  ,
\]
\emph{for any }$u\in H_{0}^{1}\left(  \Omega\right)  $\emph{\ and }$2\leq
s\leq p_{1}.$
\end{lemma}

\begin{corollary}
\noindent\ \emph{Let the assumptions of Lemma 4.1 hold. Then we have}%
\[
\left\Vert u\right\Vert _{p_{1}}^{s}\leq C\left(  ||\nabla u(t)||_{2}%
^{2}+||u(t)||_{p_{1}}^{p_{1}}\right)  ,
\]
\emph{for any }$u\in H_{0}^{1}\left(  \Omega\right)  $\emph{\ and }$2\leq
s\leq p_{1}.$
\end{corollary}

\begin{corollary}
\noindent\emph{\ Let the assumptions of Lemma 4.1 hold. Then we have}%
\[
\varrho^{\frac{s}{p_{1}}}(u)\leq C\left(  \left\vert H(t)\right\vert
+||u_{t}(t)||_{2}^{2}+\varrho(u)+\int_{\Omega}\int_{0}^{1}\int_{\tau_{1}%
}^{\tau_{2}}\tau\left(  \mu_{2}\left(  \tau\right)  +\xi(x)\right)
\frac{\left\vert z\right\vert ^{m(x)}}{m(x)}d\tau d\rho dx\right)  ,
\]
\emph{for any }$u\in H_{0}^{1}\left(  \Omega\right)  $\emph{\ and }$2\leq
s\leq p_{1}.$
\end{corollary}

\begin{corollary}
\noindent\emph{Let the assumptions of Lemma 4.1 hold. Then we have}%
\[
\left\Vert u\right\Vert _{p_{1}}^{s}\leq C\left(  \left\vert H(t)\right\vert
+||u_{t}(t)||_{2}^{2}+||u(t)||_{p_{1}}^{p_{1}}+\int_{\Omega}\int_{0}^{1}%
\int_{\tau_{1}}^{\tau_{2}}\tau\left(  \mu_{2}\left(  \tau\right)
+\xi(x)\right)  \frac{\left\vert z\right\vert ^{m(x)}}{m(x)}d\tau d\rho
dx\right)  ,
\]
\emph{for any }$u\in H_{0}^{1}\left(  \Omega\right)  $\emph{\ and }$2\leq
s\leq p_{1}.$
\end{corollary}

\begin{lemma}
\noindent\emph{Let the assumptions of Lemma 4.1 hold and let }$u$\emph{\ be
the solution of }$\left(  \ref{P1}\right)  $\emph{. Then,}
\[
\varrho(u)\geq C||u(t)||_{p_{1}}^{p_{1}}.
\]

\end{lemma}

\begin{lemma}
\emph{Let the assumptions of Lemma 4.1 hold and let }$u$\emph{\ be the
solution of }$\left(  \ref{P1}\right)  $\emph{. Then,}
\[
\int_{\Omega}\left\vert u\right\vert ^{m(x)}dx\leq C\left(  \varrho
^{\frac{m_{1}}{p_{1}}}(u)+\varrho^{\frac{m_{2}}{p_{1}}}(u)\right)  .
\]

\end{lemma}

Our main result reads as follows.

\begin{theorem}
[Main]\noindent\textbf{\ }\emph{Let the conditions }$\left(  \ref{eq1}\right)
$\emph{\ and }$\left(  \ref{eq2}\right)  $\emph{\ be fulfilled. Assume further
}$E(0)<0.$ \noindent\noindent\noindent\emph{T}\noindent\emph{hen, the solution
of }$\left(  \ref{P1}\right)  $ \emph{does not exist globally.}
\end{theorem}

\noindent\textbf{Proof.} We define
\begin{equation}
L(t):=H^{1-\alpha}(t)+\varepsilon\int_{\Omega}uu_{t}(x,t)dx,\label{eq11}%
\end{equation}
for small $\varepsilon$ to be chosen later and
\begin{equation}
0<\alpha\leq\min\left\{  \frac{p_{1}-2}{2p_{1}},\frac{p_{1}-m_{2}}%
{p_{1}\left(  m_{2}-1\right)  }\right\}  .\label{eq12}%
\end{equation}
A direct differentiation of $\left(  \ref{eq11}\right)  $ using $\left(
\ref{P1}\right)  _{1}$ gives
\begin{align*}
L^{\prime}(t)  & =\left(  1-\alpha\right)  H^{-\alpha}(t)H^{\prime
}(t)+\varepsilon\int_{\Omega}\left[  u_{t}^{2}-|\nabla u(t)|^{2}\right]
dx+\varepsilon\int_{\Omega}\left\vert u\right\vert ^{p(x)}dx\\
& -\varepsilon\mu_{1}\int_{\Omega}uu_{t}(x,t)\left\vert u_{t}(x,t)\right\vert
^{m(x)-2}dx\\
& -\varepsilon\int_{\Omega}\int_{\tau_{1}}^{\tau_{2}}\mu_{2}\left(
\tau\right)  uz(x,1,t,\tau)\left\vert z(x,1,t,\tau)\right\vert ^{m(x)-2}d\tau
dx.
\end{align*}
Using the definition of $H(t)$ and for $0<a<1,$ we have
\begin{align}
L^{\prime}(t)  & \geq C_{0}\left(  1-\alpha\right)  H^{-\alpha}(t)\left[
\int_{\Omega}\left\vert u_{t}(t)\right\vert ^{m(x)}dx+\int_{\Omega}\int%
_{\tau_{1}}^{\tau_{2}}\left\vert z(x,1,t)\right\vert ^{m(x)}d\tau dx\right]
\nonumber\\
& +\varepsilon\left(  \left(  1-a\right)  p_{1}H(t)+\frac{\left(  1-a\right)
p_{1}}{2}||u_{t}||_{2}^{2}+\frac{\left(  1-a\right)  p_{1}}{2}||\nabla
u(t)||_{2}^{2}\right) \nonumber\\
& +\varepsilon\left(  1-a\right)  p_{1}\int_{\Omega}\int_{0}^{1}\int_{\tau
_{1}}^{\tau_{2}}\tau\left(  \mu_{2}\left(  \tau\right)  +\xi(x)\right)
\frac{\left\vert z(x,\rho,t,\tau)\right\vert ^{m(x)}}{m(x)}d\tau d\rho
dx\label{eq13}\\
& +\varepsilon\int_{\Omega}\left[  u_{t}^{2}-|\nabla u(t)|_{2}^{2}\right]
dx+\varepsilon ab\int_{\Omega}\left\vert u\right\vert ^{p(x)}dx\nonumber\\
& -\varepsilon\mu_{1}\int_{\Omega}uu_{t}(x,t)\left\vert u_{t}(x,t)\right\vert
^{m(x)-2}dx\nonumber\\
& -\varepsilon\int_{\Omega}\int_{\tau_{1}}^{\tau_{2}}\mu_{2}\left(
\tau\right)  uz(x,1,t,\tau)\left\vert z(x,1,t,\tau)\right\vert ^{m(x)-2}d\tau
dx.\nonumber
\end{align}
Hence,
\begin{align}
L^{\prime}(t)  & \geq C_{0}\left(  1-\alpha\right)  H^{-\alpha}(t)\left[
\int_{\Omega}\left\vert u_{t}(t)\right\vert ^{m(x)}dx+\int_{\Omega}\int%
_{\tau_{1}}^{\tau_{2}}\left\vert z(x,1,t)\right\vert ^{m(x)}d\tau dx\right]
\nonumber\\
& +\varepsilon\left(  1-a\right)  p_{1}H(t)+\varepsilon\frac{\left(
1-a\right)  p_{1}+2}{2}||u_{t}||_{2}^{2}+\varepsilon\frac{\left(  1-a\right)
p_{1}-2}{2}||\nabla u||_{2}^{2}\label{eq14}\\
& +\varepsilon\left(  1-a\right)  p_{1}\int_{\Omega}\int_{0}^{1}\int_{\tau
_{1}}^{\tau_{2}}\tau\left(  \mu_{2}\left(  \tau\right)  +\xi(x)\right)
\frac{\left\vert z(x,\rho,t,\tau)\right\vert ^{m(x)}}{m(x)}d\tau d\rho
dx+\varepsilon a\varrho(u)\nonumber\\
& -\varepsilon\mu_{1}\int_{\Omega}uu_{t}(x,t)\left\vert u_{t}(x,t)\right\vert
^{m(x)-2}dx\nonumber\\
& -\varepsilon\int_{\Omega}\int_{\tau_{1}}^{\tau_{2}}\mu_{2}\left(
\tau\right)  uz(x,1,t,\tau)\left\vert z(x,1,t,\tau)\right\vert ^{m(x)-2}d\tau
dx.\nonumber
\end{align}
Recalling Young's inequality, we have
\begin{equation}
\int_{\Omega}\left\vert u_{t}\right\vert ^{m(x)-1}\left\vert u\right\vert
dx\leq\frac{1}{m_{1}}\int_{\Omega}\delta^{m(x)}\left\vert u\right\vert
^{m(x)}dx+\frac{m_{2}-1}{m_{2}}\int_{\Omega}\delta^{-\frac{m(x)}{m(x)-1}%
}\left\vert u_{t}\right\vert ^{m(x)}dx\label{eq15}%
\end{equation}
and
\begin{align}
& \int_{\Omega}\int_{\tau_{1}}^{\tau_{2}}\mu_{2}\left(  \tau\right)
\left\vert z(x,1,t,\tau)\right\vert ^{m(x)-1}\left\vert u\right\vert d\tau
dx\label{eq16}\\
& \leq\frac{\mu_{1}\left(  \tau_{2}-\tau_{1}\right)  }{m_{1}}\int_{\Omega
}\delta^{m(x)}\left\vert u\right\vert ^{m(x)}dx+\mu_{1}\frac{m_{2}-1}{m_{2}%
}\int_{\Omega}\int_{\tau_{1}}^{\tau_{2}}\delta^{-\frac{m(x)}{m(x)-1}%
}\left\vert z(x,1,t,\tau)\right\vert ^{m(x)}d\tau dx.\nonumber
\end{align}
As in $\cite{26}$, estimates $\left(  \ref{eq15}\right)  $ and $\left(
\ref{eq16}\right)  $ remain valid even if $\delta$ is time dependent.
Therefore, taking $\delta$ such that
\[
\delta^{-\frac{m(x)}{m(x)-1}}=kH^{-\alpha}(t),
\]
for large $k\geq1$ to be specified later, we get
\begin{align}
\int_{\Omega}\delta^{-\frac{m(x)}{m(x)-1}}\left\vert u_{t}\right\vert
^{m(x)}dx  & =kH^{-\alpha}(t)\int_{\Omega}\left\vert u_{t}\right\vert
^{m(x)}dx,\label{eq17}\\
\int_{\Omega}\int_{\tau_{1}}^{\tau_{2}}\delta^{-\frac{m(x)}{m(x)-1}}\left\vert
z(x,1,t,\tau)\right\vert ^{m(x)}d\tau dx  & =kH^{-\alpha}(t)\int_{\Omega}%
\int_{\tau_{1}}^{\tau_{2}}\left\vert z(x,1,t,\tau)\right\vert ^{m(x)}d\tau
dx\nonumber
\end{align}
and
\begin{align}
\int_{\Omega}\delta^{m(x)}\left\vert u\right\vert ^{m(x)}dx  & =\int_{\Omega
}k^{1-m(x)}H^{\alpha\left(  m(x)-1\right)  }(t)\left\vert u\right\vert
^{m(x)}dx\nonumber\\
& \leq k^{1-m_{1}}H^{\alpha\left(  m_{2}-1\right)  }(t)\int_{\Omega}\left\vert
u\right\vert ^{m(x)}dx.\label{eq18}%
\end{align}
Using Lemmas 4.5 and 4.6, we get
\[
H^{\alpha\left(  m_{2}-1\right)  }(t)\int_{\Omega}\left\vert u\right\vert
^{m(x)}dx\leq C\left[  \left(  \varrho(u)\right)  ^{\frac{m_{1}}{p_{1}}%
+\alpha\left(  m_{2}-1\right)  }+\left(  \varrho(u)\right)  ^{\frac{m_{2}%
}{p_{1}}+\alpha\left(  m_{2}-1\right)  }\right]  .
\]
By virtue of $\left(  \ref{eq12}\right)  $, we deduce that
\[
s=m_{1}+\alpha p_{1}\left(  m_{2}-1\right)  \leq p_{1}\text{ \ and \ }%
s=m_{2}+\alpha p_{1}\left(  m_{2}-1\right)  \leq p_{1}.
\]
Hence, Lemma 4.1 yields
\begin{equation}
H^{\alpha\left(  m_{2}-1\right)  }(t)\int_{\Omega}\left\vert u\right\vert
^{m(x)}dx\leq C\left(  ||\nabla u(t)||_{2}^{2}+\varrho(u)\right)
.\label{eq19}%
\end{equation}
Combining $\left(  \ref{eq14}\right)  $-$\left(  \ref{eq19}\right)  $, we
arrive at
\begin{align}
L^{\prime}(t)  & \geq\left(  1-\alpha\right)  H^{-\alpha}(t)\left[
C_{0}-\varepsilon\left(  \frac{m_{2}-1}{m_{2}}\right)  k\right]  \int_{\Omega
}\left\vert u_{t}(t)\right\vert ^{m(x)}dx\nonumber\\
& +\left(  1-\alpha\right)  H^{-\alpha}(t)\left[  C_{0}-\varepsilon\left(
\frac{m_{2}-1}{m_{2}}\right)  \mu_{1}k\right]  \int_{\Omega}\int_{\tau_{1}%
}^{\tau_{2}}\left\vert z(x,1,t,\tau)\right\vert ^{m(x)}d\tau dx\nonumber\\
& +\varepsilon\left(  \frac{\left(  p_{1}-2\right)  -ap_{1}}{2}-\frac{C}%
{m_{1}\mu_{1}\left(  \tau_{2}-\tau_{1}\right)  k^{m_{1}-1}}\right)  ||\nabla
u(t)||_{2}^{2}\label{eq20}\\
& +\varepsilon\left(  1-a\right)  p_{1}H(t)+\varepsilon\frac{\left(
1-a\right)  p_{1}+2}{2}||u_{t}||_{2}^{2}+\varepsilon\left(  a-\frac{C}%
{m_{1}\mu_{1}\left(  \tau_{2}-\tau_{1}\right)  k^{m_{1}-1}}\right)
\varrho(u)\nonumber\\
& +\varepsilon\left(  1-a\right)  p_{1}\int_{\Omega}\int_{0}^{1}\int_{\tau
_{1}}^{\tau_{2}}\tau\left(  \mu_{2}\left(  \tau\right)  +\xi(x)\right)
\frac{\left\vert z(x,\rho,t,\tau)\right\vert ^{m(x)}}{m(x)}d\tau d\rho
dx.\nonumber
\end{align}
At this point, we choose $a$ small enough so that
\[
\frac{\left(  p_{1}-2\right)  -ap_{1}}{2}>0.
\]
The constant $k$ is chosen to be large so that
\[
\frac{\left(  p_{1}-2\right)  -ap_{1}}{2}-\frac{C}{m_{1}\mu_{1}\left(
\tau_{2}-\tau_{1}\right)  k^{m_{1}-1}}>0
\]
and%
\[
a-\frac{C}{m_{1}\mu_{1}\left(  \tau_{2}-\tau_{1}\right)  k^{m_{1}-1}}>0.
\]
Once $a$ and $k$ are fixed, we pick $\varepsilon$ small enough, so that
\begin{align*}
C_{0}-\varepsilon\left(  \frac{m_{2}-1}{m_{2}}\right)  k  & >0,\\
C_{0}-\varepsilon\left(  \frac{m_{2}-1}{m_{2}}\right)  \mu_{1}k  & >0
\end{align*}
and
\[
L(0)=H^{1-\alpha}(0)+\varepsilon\int_{\Omega}u_{0}(x)u_{1}(x)dx>0.
\]
Thus, $\left(  \ref{eq20}\right)  $ takes the form
\begin{equation}
L^{\prime}(t)\geq\gamma\left[  H(t)+||u_{t}||_{2}^{2}+||\nabla u||_{2}%
^{2}+\varrho(u)+\int_{\Omega}\int_{0}^{1}\int_{\tau_{1}}^{\tau_{2}}\tau\left(
\mu_{2}\left(  \tau\right)  +\xi(x)\right)  \frac{\left\vert z\right\vert
^{m(x)}}{m(x)}d\tau d\rho dx\right]  ,\label{eq21}%
\end{equation}
for a constant $\gamma>0$. Consequently, we have
\[
L(t)\geq L(0)>0,\text{ \ \ \ }\forall t\geq0.
\]
Next, we want to show, for two constants $\sigma,\Gamma>0,$ that
\[
L^{\prime}(t)\geq\Gamma L^{\sigma}(t).
\]
For this reason, we estimate
\[
\left\vert \int_{\Omega}uu_{t}(x,t)dx\right\vert \leq\left\Vert u\right\Vert
_{2}\left\Vert u_{t}\right\Vert _{2}\leq C\left\Vert u\right\Vert _{p_{1}%
}\left\Vert u_{t}\right\Vert _{2},
\]
which implies
\[
\left\vert \int_{\Omega}uu_{t}(x,t)dx\right\vert ^{\frac{1}{1-\alpha}}\leq
C\left\Vert u\right\Vert _{p_{1}}^{\frac{1}{1-\alpha}}\left\Vert
u_{t}\right\Vert _{2}^{\frac{1}{1-\alpha}}
\]
and Young's inequality yields
\[
\left\vert \int_{\Omega}uu_{t}(x,t)dx\right\vert ^{\frac{1}{1-\alpha}}\leq
C\left[  \left\Vert u\right\Vert _{p_{1}}^{\frac{\mu}{1-\alpha}}+\left\Vert
u_{t}\right\Vert _{2}^{\frac{\theta}{1-\alpha}}\right]  ,
\]
where $1/\mu+1/\theta=1.$ The choice of $\theta=2\left(  1-\alpha\right)  $
will make $\mu/\left(  1-\alpha\right)  =2/\left(  1-2\alpha\right)  \leq
p_{1}$ by $\left(  \ref{eq12}\right)  $. Therefore,
\[
\left\vert \int_{\Omega}uu_{t}(x,t)dx\right\vert ^{\frac{1}{1-\alpha}}\leq
C\left[  \left\Vert u\right\Vert _{p_{1}}^{s}+\left\Vert u_{t}\right\Vert
_{2}^{2}\right]  ,
\]
where $s=\mu/\left(  1-\alpha\right)  .$ Using Corollary 4.4, we have
\begin{align*}
& \left\vert \int_{\Omega}uu_{t}(x,t)dx\right\vert ^{\frac{1}{1-\alpha}}\\
& \leq C\left[  H(t)+||u_{t}||_{2}^{2}+\varrho(u)+\int_{\Omega}\int_{0}%
^{1}\int_{\tau_{1}}^{\tau_{2}}\tau\left(  \mu_{2}\left(  \tau\right)
+\xi(x)\right)  \frac{\left\vert z\right\vert ^{m(x)}}{m(x)}d\tau d\rho
dx\right]  .
\end{align*}
On the other hand,
\begin{align*}
L^{1/(1-\alpha)}(t)  & =\left[  H^{\left(  1-\alpha\right)  }(t)+\varepsilon
\int_{\Omega}uu_{t}(x,t)dx\right]  ^{1/(1-\alpha)}\\
& \leq2^{1/(1-\alpha)}\left[  H(t)+\left\vert \int_{\Omega}uu_{t}dx\right\vert
^{1/(1-\alpha)}\right] \\
& \leq C\left[  H(t)+||u_{t}(t)||_{2}^{2}+\varrho(u)+\int_{\Omega}\int_{0}%
^{1}\int_{\tau_{1}}^{\tau_{2}}\tau\left(  \mu_{2}\left(  \tau\right)
+\xi(x)\right)  \frac{\left\vert z\right\vert ^{m(x)}}{m(x)}d\tau d\rho
dx\right]  .
\end{align*}
Hence, $\left(  \ref{eq21}\right)  $ implies, for some $\chi>0,$
\[
L^{\prime}(t)\geq\chi L^{1/(1-\alpha)}(t).
\]
Integration over $(0,t)$ yields
\[
L^{\alpha/(1-\alpha)}(t)\geq\frac{1}{L^{-\alpha/(1-\alpha)}(0)-\chi\alpha
t/(1-\alpha)}.
\]
Therefore, the solutions cannot be exist after the time
\[
\widetilde{T}\leq\frac{1-\alpha}{\chi\alpha\left[  L(0)\right]  ^{\alpha
/(1-\alpha)}}.
\]
This completes the proof.$\ \ \ \ \ \ \ \ \ \ \
\endproof
$

Next, we estimate the lower bound for the existence time or the life span of
the solution. For this purpose, we set\noindent%
\[
\phi(t)=\int_{\Omega}\frac{\left\vert u(x,t)\right\vert ^{p(x)}}{p(x)}dx.
\]

\begin{theorem}
\noindent\emph{Assume }$\left(  \ref{eq1}\right)  $\emph{\ holds and that }$u$
\emph{is a solution of }$\left(  \ref{P1}\right)  $\emph{\ that does not exist
after a time }$\widetilde{T}$\emph{. Then we have}
\[
\int_{\phi(0)}^{\infty}\frac{dy}{cy^{p_{2}-1}+cy^{p_{1}-1}+y+E(0)}%
\leq\widetilde{T}.
\]

\end{theorem}

\noindent\noindent\textbf{Proof.} A direct differentiation of $\phi(t)$
yields
\[
\phi^{\prime}(t)=\int_{\Omega}\left\vert u(x,t)\right\vert ^{p(x)-2}%
u(x,t)u_{t}(x,t)dx.
\]
Using Young's inequality, we have
\begin{equation}
\phi^{\prime}(t)\leq\frac{1}{2}\int_{\Omega}u_{t}^{2}dx+\frac{1}{2}%
\int_{\Omega}\left\vert u(x)\right\vert ^{2p(x)-2}dx.\label{eq22}%
\end{equation}
\noindent To estimate the second term of $\left(  \ref{eq22}\right)  $, we
define
\[
\Omega_{+}=\left\{  x\in\Omega\text{ }|\text{ }\left\vert u(x)\right\vert
\geq1\right\}  \text{ \ and \ }\Omega_{-}=\left\{  x\in\Omega\text{ }|\text{
}\left\vert u(x)\right\vert <1\right\}  .
\]
Thus, we have
\begin{align*}
\int_{\Omega}\left\vert u(x)\right\vert ^{2p(x)-2}dx  & =\int_{\Omega_{+}%
}\left\vert u(x)\right\vert ^{2p(x)-2}dx+\int_{\Omega_{-}}\left\vert
u(x)\right\vert ^{2p(x)-2}dx\\
& \leq\int_{\Omega_{+}}\left\vert u(x)\right\vert ^{2p_{2}-2}dx+\int%
_{\Omega_{-}}\left\vert u(x)\right\vert ^{2p_{1}-2}dx\\
& \leq\int_{\Omega}\left\vert u(x)\right\vert ^{2p_{2}-2}dx+\int_{\Omega
}\left\vert u(x)\right\vert ^{2p_{1}-2}dx\\
& =\left\Vert u\right\Vert _{2p_{2}-2}^{2p_{2}-2}+\left\Vert u\right\Vert
_{2p_{1}-2}^{2p_{1}-2}.
\end{align*}
Using the embedding $H_{0}^{1}(\Omega)\hookrightarrow L^{p_{1}}(\Omega)$ and
$H_{0}^{1}(\Omega)\hookrightarrow L^{p_{2}}(\Omega)$, we have
\[
\frac{1}{2}\int_{\Omega}\left\vert u(x)\right\vert ^{2p(x)-2}dx\leq c\left[
\left\Vert \nabla u\right\Vert _{2}^{2(p_{2}-1)}+\left\Vert \nabla
u\right\Vert _{2}^{2(p_{1}-1)}\right]  .
\]
Therefore, $\left(  \ref{eq22}\right)  $ becomes
\begin{align*}
\phi^{\prime}(t)  & \leq\frac{1}{2}\left\Vert u_{t}\right\Vert _{2}%
^{2}+c\left[  \left(  2\phi(t)\right)  ^{p_{2}-1}+\left(  2\phi(t)\right)
^{p_{1}-1}\right] \\
& \leq\phi(t)+c\left(  \phi(t)\right)  ^{p_{2}-1}+c\left(  \phi(t)\right)
^{p_{1}-1}+E(0).
\end{align*}
Integrating both sides of the last inequality over $\left(  0,\widetilde{T}%
\right)  ,$ we obtain
\[
\int_{\phi(0)}^{\infty}\frac{dy}{y+E(0)+cy^{p_{2}-1}+cy^{p_{1}-1}}%
\leq\widetilde{T},
\]
which is the desired result.$\ \ \ \ \ \ \ \ \ \ \
\endproof
$

\section{Global Existance and Decay}

In this section, we show that any global solution decays either exponentially
or polynomially depending on the value of $m(\cdot).$ For this purpose, we
strat by showing that the solution exists globally under some conditions on
the initial data.

Let us assume that%
\begin{equation}
I\left(  t\right)  :=||\nabla u||_{2}^{2}-\int_{\Omega}\left\vert u\right\vert
^{p(x)}dx\label{q23}%
\end{equation}
and%
\begin{align}
J\left(  t\right)   & :=\frac{1}{2}||\nabla u||_{2}^{2}+\int_{\Omega}\int%
_{0}^{1}\int_{\tau_{1}}^{\tau_{2}}\tau\left(  \mu_{2}\left(  \tau\right)
+\xi(x)\right)  \frac{\left\vert z(x,\rho,t,\tau)\right\vert ^{m(x)}}%
{m(x)}d\tau d\rho dx\label{eq24}\\
& -\int_{\Omega}\frac{\left\vert u\right\vert ^{p(x)}}{p(x)}dx.\nonumber
\end{align}
Hence,%
\begin{equation}
E\left(  t\right)  :=\frac{1}{2}||u_{t}||_{2}^{2}+J\left(  t\right)
.\label{eq25}%
\end{equation}

\begin{lemma}
\emph{Suppose that the initial data }$u_{0}$\emph{\ and }$u_{1}$%
\emph{\ satisfying }$I\left(  0\right)  >0$ \emph{and}%
\begin{equation}
\beta=c_{p}\left[  \left(  \frac{2p_{1}}{p_{1}-2}E(0)\right)  ^{\frac{p_{2}%
-2}{2}}+\left(  \frac{2p_{1}}{p_{1}-2}E(0)\right)  ^{\frac{p_{1}-2}{2}%
}\right]  <\frac{p_{1}-2}{2p_{1}}<1,\label{eq26}%
\end{equation}
\emph{where }$c_{p}$\emph{\ is the Poincar\`{e}'s constant. Then }$I\left(
t\right)  >0,$\emph{\ for }$t\in\left[  0,T\right]  .$
\end{lemma}

\noindent\textbf{Proof.} Since $I(0)>0$ we deduce by continuity that there
exists $\widetilde{T}\leq T$ such that $I(t)\geq0$ for all $t\in\left[
0,\widetilde{T}\right]  $. This implies that, for all $t\in\left[
0,\widetilde{T}\right]  $, we have%
\begin{align*}
J\left(  t\right)   & \geq\frac{1}{2}||\nabla u||_{2}^{2}+\int_{\Omega}%
\int_{0}^{1}\int_{\tau_{1}}^{\tau_{2}}\tau\left(  \mu_{2}\left(  \tau\right)
+\xi(x)\right)  \frac{\left\vert z(x,\rho,t,\tau)\right\vert ^{m(x)}}%
{m(x)}d\tau d\rho dx-\frac{1}{p_{1}}\int_{\Omega}\left\vert u\right\vert
^{p(x)}dx\\
& \geq\frac{p_{1}-2}{2p_{1}}||\nabla u||_{2}^{2}+\int_{\Omega}\int_{0}^{1}%
\int_{\tau_{1}}^{\tau_{2}}\tau\left(  \mu_{2}\left(  \tau\right)
+\xi(x)\right)  \frac{\left\vert z(x,\rho,t,\tau)\right\vert ^{m(x)}}%
{m(x)}d\tau d\rho dx+\frac{1}{p_{1}}I\left(  t\right) \\
& \geq\frac{p_{1}-2}{2p_{1}}||\nabla u||_{2}^{2}.
\end{align*}
Thus%
\begin{equation}
||\nabla u||_{2}^{2}\leq\frac{2p_{1}}{p_{1}-2}J\left(  t\right)  \leq
\frac{2p_{1}}{p_{1}-2}E\left(  t\right)  \leq\frac{2p_{1}}{p_{1}-2}E\left(
0\right)  .\label{eq27}%
\end{equation}
In the other hand, using the embedding $H_{0}^{1}(\Omega)\hookrightarrow
L^{p_{1}}(\Omega)$ and $H_{0}^{1}(\Omega)\hookrightarrow L^{p_{2}}(\Omega) $
together with Poincar\`{e}'s inequality, to get%
\begin{align}
\int_{\Omega}\left\vert u\right\vert ^{p(x)}dx  & \leq\int_{\left\vert
u\right\vert \geq1}\left\vert u\right\vert ^{p_{2}}dx+\int_{\left\vert
u\right\vert <1}\left\vert u\right\vert ^{p_{1}}dx\nonumber\\
& \leq c_{p}\left(  ||\nabla u||_{2}^{p_{2}}+||\nabla u||_{2}^{p_{1}}\right)
\nonumber\\
& =c_{p}\left(  ||\nabla u||_{2}^{p_{2}-2}+||\nabla u||_{2}^{p_{1}-2}\right)
||\nabla u||_{2}^{2}\label{eq28}\\
& \leq c_{p}\left[  \left(  \frac{2p_{1}}{p_{1}-2}E(0)\right)  ^{\frac
{p_{2}-2}{2}}+\left(  \frac{2p_{1}}{p_{1}-2}E(0)\right)  ^{\frac{p_{1}-2}{2}%
}\right]  ||\nabla u||_{2}^{2}\nonumber
\end{align}
Considering $\left(  \ref{eq26}\right)  $, we infer that $I\left(  t\right)
>0,$ for all $t\in\left[  0,\widetilde{T}\right]  $. By repeating this
procedure, $\widetilde{T}$ can be extended to $T$.

\begin{theorem}
\emph{If the initial data }$u_{0}$\emph{\ and }$u_{1}$\emph{\ satisfying the
conditions of Lemma 5.1, then the solution of }$\left(  \ref{P1}\right)
$\emph{\ is uniformly bounded and global in time.}
\end{theorem}

\begin{description}
\item[Proof.] It suffices to show that $||u_{t}||_{2}^{2}+||\nabla u||_{2}%
^{2}$ is bounded independently of $t$. Clearly,%
\[
E(0)\geq E(t)=\frac{1}{2}||u_{t}||_{2}^{2}+J\left(  t\right)  \geq\frac{1}%
{2}||u_{t}||_{2}^{2}+\frac{p_{1}-2}{2p_{1}}||\nabla u||_{2}^{2}.
\]
Thus, $||u_{t}||_{2}^{2}+||\nabla u||_{2}^{2}$ is uniformaly bounded by
$cE(0)$ independently of $t$.
\end{description}

\begin{lemma}
[\noindent Komornik $\cite{28}$ p. 103 and 124]\emph{Let }$E:\mathbb{R}%
^{+}\longrightarrow\mathbb{R}^{+}$\emph{\ be a nonincreasing function. Assume
that there exist }$\sigma>0,\omega>0$\emph{\ such that}
\[
\int_{s}^{\infty}E^{1+\sigma}(t)dt\leq\frac{1}{\omega}E^{\sigma}%
(0)E(s)=cE(s),\text{ \ \ \ }\forall s>0.
\]
\emph{Then, }$\forall t\geq0,$%
\begin{align*}
E(t)  & \leq cE(0)/(1+t)^{1/\sigma},\text{ \ \ if \ }\sigma>0,\\
E(t)  & \leq cE(0)e^{-\omega t},\text{ \ \ \ \ \ \ \ \ \ \ if \ }\sigma=0.
\end{align*}

\end{lemma}

\noindent Before we state the main theorem, we need the following technical lemma.

\begin{lemma}
\noindent\emph{The functional}
\[
F(t)=\int_{\Omega}\int_{0}^{1}\int_{\tau_{1}}^{\tau_{2}}\tau e^{-\rho\tau
}\left(  \mu_{2}\left(  \tau\right)  +\xi(x)\right)  \frac{\left\vert
z(x,\rho,t,\tau)\right\vert ^{m(x)}}{m(x)}d\tau d\rho dx,
\]
\emph{satisfies, along the solution of }$\left(  \ref{P1}\right)  $ \emph{and
for some two positive constants }$\alpha_{1},\alpha_{2},$ \emph{,}
\begin{equation}
F^{\prime}(t)\leq\alpha_{1}\int_{\Omega}|u_{t}|^{m(x)}dx-\alpha_{2}%
\int_{\Omega}\int_{0}^{1}\int_{\tau_{1}}^{\tau_{2}}\left(  \mu_{2}\left(
\tau\right)  +\xi(x)\right)  |z|^{m(x)}d\tau d\rho dx.\label{eq29}%
\end{equation}

\end{lemma}

\noindent\textbf{Proof.} A direct differentiation of $F(t),$ using $\left(
\ref{P1}\right)  _{2}$, leads to \noindent%
\begin{align*}
F^{\prime}(t)  & =-\int_{\Omega}\int_{\tau_{1}}^{\tau_{2}}\left(  \mu
_{2}\left(  \tau\right)  +\xi(x)\right)  \int_{0}^{1}e^{-\rho\tau}\left\vert
z\right\vert ^{m(x)-1}z_{\rho}d\rho d\tau dx\\
& =-\int_{\Omega}\int_{\tau_{1}}^{\tau_{2}}\frac{\mu_{2}\left(  \tau\right)
+\xi(x)}{m(x)}\int_{0}^{1}\left\{  \frac{\partial}{\partial\rho}\left[
e^{-\rho\tau}\left\vert z\right\vert ^{m(x)}\right]  +\tau e^{-\rho\tau
}\left\vert z\right\vert ^{m(x)}\right\}  d\rho d\tau dx\\
& =-\int_{\Omega}m(x)\int_{\tau_{1}}^{\tau_{2}}\frac{\mu_{2}\left(
\tau\right)  +\xi(x)}{m(x)}\left[
\begin{array}
[c]{c}%
e^{-\tau}\left\vert z(x,1,t,\tau)\right\vert ^{m(x)}-\left\vert z(x,0,t,\tau
)\right\vert ^{m(x)}\\
+\tau\int_{0}^{1}e^{-\rho\tau}\left\vert z\right\vert ^{m(x)}d\rho
\end{array}
\right]  d\tau dx\\
& \leq\alpha_{1}\int_{\Omega}|u_{t}|^{m(x)}dx-\alpha_{2}\int_{\Omega}\int%
_{0}^{1}\int_{\tau_{1}}^{\tau_{2}}\frac{\mu_{2}\left(  \tau\right)  +\xi
(x)}{m(x)}|z|^{m(x)}d\tau d\rho dx.\ \ \ \ \ \ \ \ \ \ \
\endproof
\end{align*}

Our main result reads as follows.\newline

\begin{theorem}
[Main]\emph{Assume that the conditions }$\left(  \ref{eq1}\right)
$\emph{\ and }$\left(  \ref{eq2}\right)  $\emph{\ are satisfied}. \emph{Then
there exist two positive constants }$c$ \emph{and }$\alpha$\emph{\ such that
any global solution of }$\left(  \ref{P1}\right)  $\emph{\ satisfies}%
\[
\left\{
\begin{tabular}
[c]{lll}%
$E(t)\leq cE(0)/(1+t)^{2/(m_{2}-2)},$ & if & $m_{2}>2,$\\
$E(t)\leq ce^{-\alpha t},$ & if & $m(\cdot)=2.$%
\end{tabular}
\right.
\]

\end{theorem}

\noindent\noindent\textbf{Proof.} Multiply $\left(  \ref{P1}\right)  _{1}$ by
$uE^{q}(t),$ for $q>0$ to be specified later, and integrate over $\Omega
\times\left(  s,T\right)  ,$ $s<T,$ to obtain%
\begin{align*}
& \int_{s}^{T}E^{q}(t)\int_{\Omega}\left(  uu_{tt}-u\Delta u+\mu_{1}%
uu_{t}\left\vert u_{t}\right\vert ^{m(x)-2}\right. \\
& \left.  +u\int_{\tau_{1}}^{\tau_{2}}\mu_{2}\left(  \tau\right)
z(1)\left\vert z(1)\right\vert ^{m(x)-2}d\tau-u\left\vert u\right\vert
^{p(x)-2}\right)  dxdt=0.
\end{align*}
which gives
\begin{align}
& \int_{s}^{T}E^{q}(t)\int_{\Omega}\left(  \frac{d}{dt}\left(  uu_{t}\right)
-u_{t}^{2}+|\nabla u|^{2}+\mu_{1}uu_{t}\left\vert u_{t}\right\vert
^{m(x)-2}\right. \nonumber\\
& \left.  +u\int_{\tau_{1}}^{\tau_{2}}\mu_{2}\left(  \tau\right)
z(1)\left\vert z(1)\right\vert ^{m(x)-2}d\tau-u\left\vert u\right\vert
^{p(x)-2}\right)  dxdt=0.\label{eq30}%
\end{align}
Recalling the definition of $E(t)$ given in $\left(  \ref{eq4}\right)  $,
adding and subtracting some term and using the relation
\[
\frac{d}{dt}\left(  E^{q}(t)\int_{\Omega}uu_{t}dx\right)  =qE^{q-1}%
(t)E^{\prime}(t)\int_{\Omega}uu_{t}dx+E^{q}(t)\frac{d}{dt}\int_{\Omega}%
uu_{t}dx,
\]
Equation $\left(  \ref{eq30}\right)  $ becomes
\begin{align}
& 2\int_{s}^{T}E^{q+1}(t)dt\nonumber\\
& =-\int_{s}^{T}\frac{d}{dt}\left(  E^{q}(t)\int_{\Omega}uu_{t}dx\right)
+q\int_{s}^{T}E^{q-1}(t)E^{\prime}(t)\int_{\Omega}uu_{t}dxdt\nonumber\\
& +2\int_{s}^{T}E^{q}(t)\int_{\Omega}u_{t}^{2}dx-\mu_{1}\int_{s}^{T}%
E^{q}(t)\int_{\Omega}uu_{t}\left\vert u_{t}\right\vert ^{m(x)-2}%
dxdt\nonumber\\
& -\int_{s}^{T}E^{q}(t)\int_{\Omega}u\int_{\tau_{1}}^{\tau_{2}}\mu_{2}\left(
\tau\right)  z(1)\left\vert z(1)\right\vert ^{m(x)-2}d\tau dxdt\label{eq31}\\
& +2\int_{s}^{T}E^{q}(t)\int_{0}^{1}\int_{\Omega}\int_{\tau_{1}}^{\tau_{2}%
}\frac{\tau\left(  \mu_{2}\left(  \tau\right)  +\xi(x)\right)  }%
{m(x)}\left\vert z(x,\rho,t,\tau)\right\vert ^{m(x)}d\tau dxd\rho
dt\nonumber\\
& +\int_{s}^{T}E^{q}(t)\int_{\Omega}\left\vert u\right\vert ^{p(x)}%
dxdt-2\int_{s}^{T}E^{q}(t)\int_{\Omega}\frac{\left\vert u\right\vert ^{p(x)}%
}{p(x)}dxdt.\nonumber
\end{align}
The first term in the right hand side of $\left(  \ref{eq31}\right)  $ is
estimated as follows.
\begin{align*}
\left\vert -\int_{s}^{T}\frac{d}{dt}\left(  E^{q}(t)\int_{\Omega}%
uu_{t}dx\right)  \right\vert  & =\left\vert E^{q}(s)\int_{\Omega}%
uu_{t}(x,s)dx-E^{q}(T)\int_{\Omega}uu_{t}(x,T)dx\right\vert \\
& \leq\frac{1}{2}E^{q}(s)\left[  \int_{\Omega}u^{2}(x,s)dx+\int_{\Omega}%
u_{t}^{2}(x,s)dx\right] \\
& +\frac{1}{2}E^{q}(T)\left[  \int_{\Omega}u^{2}(x,T)dx+\int_{\Omega}u_{t}%
^{2}(x,T)dx\right] \\
& \leq\frac{1}{2}E^{q}(s)\left[  C_{p}\left\Vert \nabla u(s)\right\Vert
_{2}^{2}+2E(s)\right] \\
& +\frac{1}{2}E^{q}(T)\left[  C_{p}\left\Vert \nabla u(T)\right\Vert _{2}%
^{2}+2E(T)\right] \\
& \leq E^{q}(s)\left[  C_{P}E(s)+E(s)\right]  +E^{q}(T)\left[  C_{P}%
E(T)+E(T)\right]  ,
\end{align*}
where $C_{P}$ is the Poincar\'{e} constant. Using the fact that $E(t)$ is
decreasing, we deduce that
\begin{equation}
\left\vert -\int_{s}^{T}\frac{d}{dt}\left(  E^{q}(t)\int_{\Omega}%
uu_{t}dx\right)  \right\vert \leq cE^{q+1}(s)\leq cE^{q}(0)E(s)\leq
cE(s).\label{eq33}%
\end{equation}
Similarly, we treat the second term:
\begin{align}
\left\vert q\int_{s}^{T}E^{q-1}(t)E^{\prime}(t)\int_{\Omega}uu_{t}%
dxdt\right\vert  & \leq-q\int_{s}^{T}E^{q-1}(t)E^{\prime}(t)\left[
C_{p}E(t)+2E(t)\right] \label{eq34}\\
& \leq-c\int_{s}^{T}E^{q}(t)E^{\prime}(t)dt\leq cE^{q+1}(s)\leq
cE(s).\nonumber
\end{align}
To handle the third term, we set
\[
\Omega_{+}=\left\{  x\in\Omega\text{ }|\text{ }\left\vert u_{t}%
(x,t)\right\vert \geq1\right\}  \text{ \ and \ }\Omega_{-}=\left\{  x\in
\Omega\text{ }|\text{ }\left\vert u_{t}(x,t)\right\vert <1\right\}
\]
and use H\"{o}lder's and Young's inequalities, to get
\begin{align*}
\left\vert \int_{s}^{T}E^{q}(t)\int_{\Omega}u_{t}^{2}dx\right\vert  &
=\left\vert \int_{s}^{T}E^{q}(t)\left[  \int_{\Omega_{+}}u_{t}^{2}%
dx+\int_{\Omega_{-}}u_{t}^{2}dx\right]  \right\vert \\
& \leq c\int_{s}^{T}E^{q}(t)\left[  \left(  \int_{\Omega_{+}}\left\vert
u_{t}\right\vert ^{m_{1}}dx\right)  ^{2/m_{1}}+\left(  \int_{\Omega_{-}%
}\left\vert u_{t}\right\vert ^{m_{2}}dx\right)  ^{2/m_{2}}\right] \\
& \leq c\int_{s}^{T}E^{q}(t)\left[  \left(  \int_{\Omega}\left\vert
u_{t}\right\vert ^{m(x)}dx\right)  ^{2/m_{1}}+\left(  \int_{\Omega}\left\vert
u_{t}\right\vert ^{m(x)}dx\right)  ^{2/m_{2}}\right] \\
& \leq c\int_{s}^{T}E^{q}(t)\left[  \left(  -E^{\prime}(t)\right)  ^{2/m_{1}%
}+\left(  -E^{\prime}(t)\right)  ^{2/m_{2}}\right] \\
& \leq c\varepsilon\int_{s}^{T}\left[  E(t)\right]  ^{qm_{1}/\left(
m_{1}-2\right)  }dt+c(\varepsilon)\int_{s}^{T}\left(  -E^{\prime}(t)\right)
dt\\
& +c\varepsilon\int_{s}^{T}E^{q+1}(t)dt+c(\varepsilon)\int_{s}^{T}\left(
-E^{\prime}(t)\right)  ^{2(q+1)/m_{2}}dt.
\end{align*}
For $m_{1}>2,$ the choice of $q=\frac{m_{2}}{2}-1$ will make $\frac{qm_{1}%
}{m_{1}-2}=q+1+\frac{m_{2}-m_{1}}{m_{1}-2}.$ Hence,
\begin{align}
& \left\vert \int_{s}^{T}E^{q}(t)\int_{\Omega}u_{t}^{2}dx\right\vert
\nonumber\\
& \leq c\varepsilon\int_{s}^{T}E^{q+1}(t)dt+c\varepsilon\left[  E(0)\right]
^{\frac{m_{2}-m_{1}}{m_{1}-2}}\int_{s}^{T}\left[  E(t)\right]  ^{q+1}%
dt+c(\varepsilon)E(s)\nonumber\\
& \leq c\varepsilon\int_{s}^{T}E^{q+1}(t)dt+c(\varepsilon)E(s).\label{eq35}%
\end{align}
For the case $m_{1}=2,$ the choice of $q=\frac{m_{2}}{2}-1,$ will give a
similar result.

For the fourth term, we use Young's inequality. So, for a.e $x\in\Omega,$ we
have
\begin{align*}
& \left\vert -\mu_{1}\int_{s}^{T}E^{q}(t)\int_{\Omega}u\left\vert
u_{t}\right\vert ^{m(x)-1}dxdt\right\vert \\
& \leq\varepsilon\int_{s}^{T}E^{q}(t)\int_{\Omega}\left\vert u(t)\right\vert
^{m(x)}dxdt+c\int_{s}^{T}E^{q}(t)\int_{\Omega}c_{\varepsilon}(x)\left\vert
u_{t}(t)\right\vert ^{m(x)}dxdt\\
& \leq\varepsilon\int_{s}^{T}E^{q}(t)\left[  \int_{\Omega_{+}}\left\vert
u(t)\right\vert ^{m_{1}}dxdt+\int_{\Omega_{-}}\left\vert u(t)\right\vert
^{m_{2}}dxdt\right] \\
& +c\int_{s}^{T}E^{q}(t)\int_{\Omega}c_{\varepsilon}(x)\left\vert
u_{t}(t)\right\vert ^{m(x)}dxdt,
\end{align*}
where we used Young's inequality with
\[
p(x)=\frac{m(x)}{m(x)-1}\text{ \ \ and \ \ }p^{\prime}(x)=m(x)
\]
and, hence,
\[
c_{\varepsilon}(x)=\varepsilon^{1-m(x)}\left(  m(x)^{-m(x)}\left(
m(x)-1\right)  \right)  ^{m(x)-1}.
\]
Therefore, using the embedding of $H_{0}^{1}(\Omega)\hookrightarrow L^{m_{1}%
}(\Omega)$ and $H_{0}^{1}(\Omega)\hookrightarrow L^{m_{2}}(\Omega) $, we
arrive at
\begin{align}
& \left\vert -\mu_{1}\int_{s}^{T}E^{q}(t)\int_{\Omega}u\left\vert
u_{t}\right\vert ^{m(x)-1}dxdt\right\vert \nonumber\\
& \leq\varepsilon\int_{s}^{T}E^{q}(t)\left[  c\left\Vert \nabla
u(s)\right\Vert _{2}^{m_{1}}+c\left\Vert \nabla u(s)\right\Vert _{2}^{m_{2}%
}\right] \nonumber\\
& +c\int_{s}^{T}E^{q}(t)\int_{\Omega}c_{\varepsilon}(x)\left\vert
u_{t}(t)\right\vert ^{m(x)}dxdt\label{eq36}\\
& \leq\varepsilon\int_{s}^{T}E^{q}(t)\left[  cE^{\frac{m_{1}-2}{2}%
}(0)E(t)+cE^{\frac{m_{2}-2}{2}}(0)E(t)\right] \nonumber\\
& +c\int_{s}^{T}E^{q}(t)\int_{\Omega}c_{\varepsilon}(x)\left\vert
u_{t}(t)\right\vert ^{m(x)}dxdt\nonumber\\
& \leq c\varepsilon\int_{s}^{T}E^{q+1}(t)+\int_{s}^{T}E^{q}(t)\int_{\Omega
}c_{\varepsilon}(x)\left\vert u_{t}(t)\right\vert ^{m(x)}dxdt.\nonumber
\end{align}
The fifth term of $\left(  \ref{eq31}\right)  $ can be estimated in a similar
manner to reach
\begin{align}
& \left\vert \int_{s}^{T}E^{q}(t)\int_{\Omega}u\int_{\tau_{1}}^{\tau_{2}}%
\mu_{2}\left(  \tau\right)  \left\vert z(x,1,t,\tau)\right\vert ^{m(x)-1}d\tau
dxdt\right\vert \nonumber\\
& \leq\varepsilon\int_{s}^{T}E^{q}(t)\left[  c\left\Vert \nabla
u(s)\right\Vert _{2}^{m_{1}}+c\left\Vert \nabla u(s)\right\Vert _{2}^{m_{2}%
}\right] \label{eq37}\\
& +c\int_{s}^{T}E^{q}(t)\int_{\Omega}c_{\varepsilon}(x)\int_{\tau_{1}}%
^{\tau_{2}}\left\vert z(x,1,t,\tau)\right\vert ^{m(x)}d\tau dxdt\nonumber\\
& \leq c\varepsilon\int_{s}^{T}E^{q+1}(t)dt+\int_{s}^{T}E^{q}(t)\int_{\Omega
}c_{\varepsilon}(x)\int_{\tau_{1}}^{\tau_{2}}\left\vert z(x,1,t,\tau
)\right\vert ^{m(x)}d\tau dxdt.\nonumber
\end{align}
The sixth term can be estimated, using Lemma 5.4, as follows,
\begin{align*}
& 2\int_{s}^{T}E^{q}(t)\int_{0}^{1}\int_{\Omega}\int_{\tau_{1}}^{\tau_{2}}%
\tau\left(  \mu_{2}\left(  \tau\right)  +\xi(x)\right)  \frac{\left\vert
z(x,\rho,t,\tau)\right\vert ^{m(x)}}{m(x)}d\tau dxd\rho dt\\
& \leq\frac{2\tau_{2}}{m_{1}}\int_{s}^{T}E^{q}(t)\int_{0}^{1}\int_{\Omega}%
\int_{\tau_{1}}^{\tau_{2}}\left(  \mu_{2}\left(  \tau\right)  +\xi(x)\right)
\left\vert z(x,\rho,t,\tau)\right\vert ^{m(x)}d\tau dxd\rho dt\\
& \leq-\frac{2\tau_{2}}{m_{1}\alpha_{2}}\int_{s}^{T}E^{q}(t)\frac{d}%
{dt}\left(  \int_{0}^{1}\int_{\Omega}\int_{\tau_{1}}^{\tau_{2}}e^{-\rho\tau
}\tau\left(  \mu_{2}\left(  \tau\right)  +\xi(x)\right)  |z|^{m(x)}d\tau
dxd\rho\right) \\
& +\frac{2\tau_{2}\alpha_{1}}{m_{1}\alpha_{2}}\int_{s}^{T}E^{q}(t)\int%
_{\Omega}|u_{t}|^{m(x)}dxdt\\
& \leq-\frac{2\tau_{2}}{m_{1}\alpha_{2}}\left[  E^{q}(t)\int_{0}^{1}%
\int_{\Omega}\int_{\tau_{1}}^{\tau_{2}}e^{-\rho\tau}\tau\left(  \mu_{2}\left(
\tau\right)  +\xi(x)\right)  |z|^{m(x)}d\tau dxd\rho\right]  _{t=s}^{t=T}\\
& +\frac{2\tau_{2}\alpha_{1}}{m_{1}\alpha_{2}}\int_{s}^{T}E^{q}(t)\int%
_{\Omega}|u_{t}|^{m(x)}dxdt.
\end{align*}
As $\xi\left(  x\right)  ,\mu_{2}\left(  \tau\right)  ,\tau$ and $\rho$ are
bounded and using $\left(  \ref{eq4}\right)  $, we obtain, for $c>0,$%
\begin{align}
& 2\int_{s}^{T}E^{q}(t)\int_{0}^{1}\int_{\Omega}\int_{\tau_{1}}^{\tau_{2}}%
\tau\left(  \mu_{2}\left(  \tau\right)  +\xi(x)\right)  \frac{\left\vert
z(x,\rho,t,\tau)\right\vert ^{m(x)}}{m(x)}d\tau dxd\rho dt\nonumber\\
& \leq cE^{q}(s)E(s)+cE^{q+1}(T)\label{eq38}\\
& \leq cE^{q}(0)E(s)+cE^{q}(T)E(s)\nonumber\\
& \leq cE(s).\nonumber
\end{align}
With the aid of $\left(  \ref{eq26}\right)  $, we can handle the last two
terms as follows.%
\begin{align}
& \int_{s}^{T}E^{q}(t)\int_{\Omega}\left\vert u\right\vert ^{p(x)}%
dxdt-2\int_{s}^{T}E^{q}(t)\int_{\Omega}\frac{\left\vert u\right\vert ^{p(x)}%
}{p(x)}dxdt\nonumber\\
& \leq\int_{s}^{T}E^{q}(t)\int_{\Omega}\left\vert u\right\vert ^{p(x)}%
dxdt-\frac{2}{p_{2}}\int_{s}^{T}E^{q}(t)\int_{\Omega}\left\vert u\right\vert
^{p(x)}dxdt\nonumber\\
& =\left(  \frac{p_{2}-2}{p_{2}}\right)  \int_{s}^{T}E^{q}(t)\int_{\Omega
}\left\vert u\right\vert ^{p(x)}dxdt\nonumber\\
& \leq\left(  \frac{p_{2}-2}{p_{2}}\right)  \beta\int_{s}^{T}E^{q}%
(t)\left\Vert \nabla u(t)\right\Vert _{2}^{2}dt\label{eq39}\\
& \leq\left(  \frac{p_{2}-2}{p_{2}}\right)  \left(  \frac{2p_{1}}{p_{1}%
-2}\right)  \beta\int_{s}^{T}E^{q+1}(t)dt\nonumber\\
& \leq\left(  \frac{p_{2}-2}{p_{2}}\right)  \int_{s}^{T}E^{q+1}(t)dt\nonumber
\end{align}
Combining $\left(  \ref{eq31}\right)  -\left(  \ref{eq39}\right)  $ and for
some $c>0$, we arrive at
\begin{align}
& \left(  2-\frac{p_{2}-2}{p_{2}}-\varepsilon\right)  \int_{s}^{T}%
E^{q+1}(t)dt\label{eq40}\\
& \leq cE(s)+\int_{s}^{T}E^{q}(t)\int_{\Omega}c_{\varepsilon}(x)\left\vert
u_{t}(t)\right\vert ^{m(x)}dxdt\nonumber\\
& +\int_{s}^{T}E^{q}(t)\int_{\Omega}c_{\varepsilon}(x)\int_{\tau_{1}}%
^{\tau_{2}}\left\vert z(x,1,t,\tau)\right\vert ^{m(x)}d\tau dxdt.\nonumber
\end{align}
As $\frac{p_{2}-2}{p_{2}}<1,$then the choice of $\varepsilon$ small enough
gives
\begin{align}
\int_{s}^{T}E^{q+1}(t)dt  & \leq cE(s)+c\int_{s}^{T}E^{q}(t)\int_{\Omega
}c_{\varepsilon}(x)\left\vert u_{t}(t)\right\vert ^{m(x)}dxdt\label{eq41}\\
& +\int_{s}^{T}E^{q}(t)\int_{\Omega}c_{\varepsilon}(x)\int_{\tau_{1}}%
^{\tau_{2}}\left\vert z(x,1,t,\tau)\right\vert ^{m(x)}d\tau dxdt.\nonumber
\end{align}
Once $\varepsilon$ is fixed, then $c_{\varepsilon}(x)$ becomes bounded (i.e.
$c_{\varepsilon}(x)\leq M$) as $m(x)$ is bounded. So, $\left(  \ref{eq41}%
\right)  $ yields, for $c>0,$
\begin{align}
\int_{s}^{T}E^{q+1}(t)dt  & \leq cE(s)+M\int_{s}^{T}E^{q}(t)\int_{\Omega
}\left\vert u_{t}(t)\right\vert ^{m(x)}dxdt\nonumber\\
& +M\int_{s}^{T}E^{q}(t)\int_{\Omega}\int_{\tau_{1}}^{\tau_{2}}\left\vert
z(x,1,t,\tau)\right\vert ^{m(x)}d\tau dxdt.\nonumber\\
& \leq cE(s)-C_{0}M\int_{s}^{T}E^{q}(t)E^{\prime}(t)dt\label{eq42}\\
& \leq cE(s)+\frac{C_{0}M}{q+1}\left[  E^{q+1}(s)-E^{q+1}(T)\right]  \leq
cE(s).\nonumber
\end{align}
\noindent As $T\longrightarrow\infty$, from $\left(  \ref{eq42}\right)  $, we
get
\[
\int_{s}^{\infty}E^{q+1}(t)dt\leq cE(s).
\]
Therefore, Komornik's Lemma is satisfied with $\sigma=q=\frac{m_{2}}{2}-1$
which implies the desired result.$\ \ \ \ \ \ \ \ \ \ \
\endproof
$

\begin{description}
\item[\textbf{Acknowledgment:}] The author would like to express his sincere
thanks to King Fahd University of Petroleum and Minerals and the
Interdisciplinary Research Center in Construction and Building Materials for
their support.
\end{description}

\end{document}